\documentstyle[12pt]{article}
\date{}
\author{Valerii Dryuma\thanks{Work supported in part by Grants RFFI and DAAD-2013}\\[5mm]
{\it Institute of Mathematics and Informatics, AS RM,}\\[3mm] {\it
5 Academiei Street, 2028 Kishinev, Moldova},\\[3mm]{\it e-mail:
valdryum@gmail.com; valerii.dryuma22@mail.ru} }
\title{ Six-dimensional spaces defined\\[1mm] by the equations KN and KdV}
\newtheorem{theorem}{Theorem}
\newtheorem{pr}{Proposition}
\newtheorem{rem}{Remark}
\begin{document}
\maketitle
\date{}
\maketitle
\begin{abstract}
\ \ \ \ On base of  three-dimensional flat metrics  obtained  with the help of solutions of  the KdV-equation
were constructed the  examples of  six-dimensional  metrics,  which are  determined by the help  of solutions
 of the Krichever-Novikov equation and the  KdV  equation. Their properties are discussed.
\end{abstract}


\medskip
\section{Flat $6D$-spaces}
\begin{theorem}

      The metric of the $6D$-space in local coordinates $x^i=(x,y,z,u,v,w)$
      \begin{equation} \label{dr:eq01}
{{\it ds}}^{2}={y}^{2}{{\it dx}}^{2}+2\,\left ({y}^{2}l(x,z)-1/2
\right ){\it dx}\,{\it dz}+2\,{\it dy}\,{\it dz}+\left ({y}^{2}\left (
l(x,z)\right )^{2}-2\,y{\frac {\partial }{\partial x}}l(x,z)+l(x,z)
\right ){{\it dz}}^{2}+$$$$+2\,{\frac {B(x,y,z){\it du}\,{\it dv}}{\left (w
u+vw+uv\right )^{2}}}+2\,{\frac {B(x,y,z){\it du}\,{\it dw}}{\left (wu
+vw+uv\right )^{2}}}+2\,{\frac {B(x,y,z){\it dv}\,{\it dw}}{\left (wu+
vw+uv\right )^{2}}}
 \end{equation}
is a flat if the following conditions on the coefficients hold
   \begin{equation}\label{dr:eq02}
B(x,y,z)=1/4\,\left ({\it F_1}(x,z)\right )^{2}{y}^{2}+1/2\,{\it F_1
}(x,z)y{\it F_2}(x,z)+1/4\,\left ({\it F_2}(x,z)\right )^{2}$$$$
{\it F_1}(x,z)=-2\,{\frac {\partial }{\partial x}}{\it F_2}(x,z),$$$$
l(x,z)=-1/3\,{\frac {{\frac {\partial }{\partial z}}{\it F_2}(x,z)-2
\,{\frac {\partial ^{3}}{\partial {x}^{3}}}{\it F_2}(x,z)}{{\frac {
\partial }{\partial x}}{\it F_2}(x,z)}},
\end{equation}
where the function ${\it F_2}(x,z)$ satisfies the equation
\begin{equation}\label{krnov1}
{\frac {\partial }{\partial z}}{\it F_2}(x,z)=-{\frac {\partial ^{3}}
{\partial {x}^{3}}}{\it F_2}(x,z)+3/2\,{\frac {\left ({\frac {
\partial ^{2}}{\partial {x}^{2}}}{\it F_2}(x,z)\right )^{2}}{{\frac {
\partial }{\partial x}}{\it F_2}(x,z)}},
\end{equation}
which is the particular case of the Krichever-Novikov equation integrable by the $IST$-method [1].
\end{theorem}

 This theorem is generalization of the following result obtained by the author [2]
\begin{theorem}
     The metric of the $3D$-space in local coordinates $x^i=(x,y,z)$
\begin{equation}\label{metr3d}
{{\it ds}}^{2}={y}^{2}{{\it dx}}^{2}+2\,\left ({y}^{2}l(x,z)-1/2
\right ){\it dx}\,{\it dz}+2\,{\it dy}\,{\it dz}+\left ({y}^{2}\left (
l(x,z)\right )^{2}-2\,y{\frac {\partial }{\partial x}}l(x,z)+l(x,z)
\right ){{\it dz}}^{2}
\end{equation}
is a flat if the function $l(x,z)$ satisfies the $KdF$-equation
\begin{equation}\label{kdf2}
{\frac {\partial }{
\partial z}}l(x,z)-3\,l(x,z){\frac {\partial }{\partial x}}l(x,z)+{\frac {\partial ^{3}}{\partial {x}^{3}}}l(x,z)=0.
\end{equation}
\end{theorem}

   Both equations can be of used to study the properties of the six-dimensional metrics  having applications in the theory
   of Calaby-Yau spaces.
\begin{theorem}
  Six-dimensional space in local coordinates $x^i=(x,y,z,u,v,w)$ equipped with the metric
\begin{equation}\label{calyau}
ds^{2}={y}^{2}{d{{x}}}^{2}+2\,\left ({y}^{2}l(x,z)+m(x,z)\right )d{{x}}d{{
z}}+$$$$+2\,d{{y}}d{{z}}+\left ({y}^{2}\left (l(x,z)\right )^{2}-2\,y{
\frac {\partial }{\partial x}}l(x,z)+2\,l(x,z)+2\,l(x,z)m(x,z)\right )
{d{{z}}}^{2}+$$$$+{v}^{2}{d{{u}}}^{2}+2\,\left ({v}^{2}L(u,w)+M(u,w)
\right )d{{u}}d{{w}}+2\,d{{v}}d{{w}}+$$$$+\left({v}^{2}\left (L(u,w)
\right )^{2}-2\,v{\frac {\partial }{\partial u}}L(u,w)+2\,L(u,w)+2\,L(
u,w)M(u,w)\right ){d{{w}}}^{2}
\end{equation}
with the coefficients
$$
l(x,z)=\left(-1/2\,{\frac {\left ({\frac {\partial ^{2}}{\partial {x}^
{2}}}{\it F_2}(x,z)\right )^{2}}{{\frac {\partial }{\partial x}}{\it
F_2}(x,z)}}+{\frac {\partial ^{3}}{\partial {x}^{3}}}{\it F_2}(x,z)
\right)\left ({\frac {\partial }{\partial x}}{\it F_2}(x,z)\right )^{
-1},$$$$
m(x,z)=-1/2,~
{\it F_1}(x,z)=-2\,{\frac {\partial }{\partial x}}{\it F_2}(x,z)
$$
and with condition
$$
{\frac {\partial }{\partial w}}M(u,w)=L(u,w){\frac {\partial }{
\partial u}}M(u,w)+{\frac {\partial }{\partial u}}L(u,w)+2\,\left ({
\frac {\partial }{\partial u}}L(u,w)\right )M(u,w)
$$
is a flat if the function $F_2(x,z)$ satisfies the KN-equation
$$
{\frac {\partial }{\partial z}}{\it F_2}(x,z)+{\frac {\partial ^{3}}
{\partial {x}^{3}}}{\it F_2}(x,z)-3/2\,{\frac {\left ({\frac {
\partial ^{2}}{\partial {x}^{2}}}{\it F_2}(x,z)\right )^{2}}{{\frac {
\partial }{\partial x}}{\it F_2}(x,z)}}=0,
$$
and the function $L(u,w)$ satisfies  the KdV-equation
$$
{\frac {\partial }{\partial w}}L(u,w)-3\,L(u,w){\frac {\partial }{
\partial u}}L(u,w)+{\frac {\partial ^{3}}{\partial {u}^{3}}}L(u,w)=0.
$$
\end{theorem}

\begin{rem}

    According  results of the authors [3]  solutions of the KN-equation can be produced by the formula
$$
(F_{2j+1})_{x}=\frac{ (F_{2j})^{2}}{(F_{2j})_x},~~j= 0,1,2,...,
$$
in which the function $F_{20}(x,z)$ is the  some initial solution of the equation (\ref{krnov1}).

   For example for the simplest function   $
    F_{20}=x$  we  get the solution
$$
 F_{21}(x,z)=1/3\,{x}^{3}+4\,z.
$$

      On the next step we find that the function
$$
F_{22}(x,z)=1/45\,{x}^{5}+4/3\,{x}^{2}z-16\,{\frac {{z}^{2}}{x}}
$$
is the solution of the $KN$-equation and so on.

        In result  we obtain sets of the six-dimensional flat spaces which are determined by solutions  of the equations KdF and KN. Some of them can be compacts and from this point of view they may be of used to further applications.
\end{rem}

\section{Ricci-flat curved 6D-space}

     Give an example of a 6-dimensional Ricci-flat space $R_{ik}=0$ with nonzero Riemann tensor $R_{i j k l}\neq 0$.

     With this aim we use the construction of Riemann extensions of the space with connection [4].

   For  the $E^n$-dimensional space in local  coordinates $x^j=(x^1,x^2,...x^n)$  with symmetric coefficients of connexion $\Pi^k_{i j}$ can be constructed  space of the Riemann
   double dimensionality $D=2n$ having  the metrics of the form
\begin{equation} \label{metr}
ds^2=-2\Pi^k_{ij} \xi^k dx^i dx^k+2\xi_k dx^k,
\end{equation}
where $\xi_k$-are additive coordinates.

    We  use this fact for construction of the $6D$- space with the Ricci- flat metric on solutions of the KdF-equation  but nonzero curvature.

    Coefficients of connection of the of the space $N^3$ with the metric (\ref{metr3d}) have the form
$$
{\it \Pi^1_{11}}=1/2\,{\frac {2\,l(x,z){y}^{2}-1}{y}},~
{\it \Pi^2_{11}}=-1/4\,{\frac {4\,{y}^{3}{\frac {\partial }{\partial x}}l(x,
z)-8\,l(x,z){y}^{2}+1}{y}},$$$$
{\it \Pi^3_{11}}=-y,~~
{\it \Pi^1_{12}}={y}^{-1},~~
{\it \Pi^2_{12}}=1/2\,{y}^{-1}
{\it \Pi^1_{13}}=1/2\,{\frac {\left (2\,l(x,z){y}^{2}-1\right )l(x,z)}{y}},$$$$
{\it \Pi^2_{13}}=-1/4\,{\frac {4\,{y}^{3}l(x,z){\frac {\partial }{\partial x
}}l(x,z)-8\,{y}^{2}\left (l(x,z)\right )^{2}+l(x,z)+4\,{y}^{2}{\frac {
\partial ^{2}}{\partial {x}^{2}}}l(x,z)-2\,y{\frac {\partial }{
\partial x}}l(x,z)}{y}},$$$$
{\it \Pi^3_{13}}=-l(x,z)y,~~
{\it \Pi^1_{23}}={\frac {l(x,z)}{y}},~~
{\it \Pi^2_{23}}=-1/2\,{\frac {2\,y{\frac {\partial }{\partial x}}l(x,z)-l(x
,z)}{y}}$$$$
{\it\Pi^1_{33}}=-1/2\,{\frac {-2\,\left ({\frac {\partial }{\partial z}}l(x
,z)\right )y+4\,yl(x,z){\frac {\partial }{\partial x}}l(x,z)-2\,{
\frac {\partial ^{2}}{\partial {x}^{2}}}l(x,z)-2\,\left (l(x,z)\right
)^{3}{y}^{2}+\left (l(x,z)\right )^{2}}{y}},$$$$
4y{\it \Pi^2_{33}}=-4\,\left (l(x,z)\right )^{2}{y}^{3}{\frac {
\partial }{\partial x}}l(x,z)-4\,l(x,z){y}^{2}{\frac {\partial ^{2}}{
\partial {x}^{2}}}l(x,z)+4\,\left ({\frac {\partial }{\partial z}}l(x,
z)\right )y+2\,{\frac {\partial ^{2}}{\partial {x}^{2}}}l(x,z)+$$$$+8\,{y}^
{2}\left ({\frac {\partial }{\partial x}}l(x,z)\right )^{2}+8\,\left (
l(x,z)\right )^{3}{y}^{2}-8\,yl(x,z){\frac {\partial }{\partial x}}l(x
,z)-\left (l(x,z)\right )^{2}-4\,{y}^{2}{\frac {\partial ^{2}}{
\partial x\partial z}}l(x,z),$$$$
{\it \Pi^3_{33}}=-\left (l(x,z)\right )^{2}y+{\frac {\partial }{\partial x}}
l(x,z).
$$

   The metrics of $6D$-space in local coordinates $x^i=(x,y,z,u,v,w)$ in considered  case is
   of the form
\begin{equation} \label{6dflat}
  {{\it ds}}^{2}=-2\,{\it \Pi^1_{11}}\,u{d{{x}}}^{2}-2\,{\it \Pi^1_{33}}\,u{d{{z}}
}^{2}-2\,{\it \Pi^3_{11}}\,w{d{{x}}}^{2}-2\,{\it \Pi^2_{11}}\,v{d{{x}}}^{2}-2\,{
\it \Pi^3_{33}}\,w{d{{z}}}^{2}-2\,{\it \Pi^2_{33}}\,v{d{{z}}}^{2}-$$$$-4\,{\it \Pi^3_{13}}
\,wd{{x}}d{{z}}-4\,{\it \Pi^1_{12}}\,ud{{x}}d{{y}}-4\,{\it \Pi^2_{12}}\,vd{{x
}}d{{y}}-4\,{\it \Pi^1_{13}}\,ud{{x}}d{{z}}-4\,{\it \Pi^3_{12}}\,vd{{x}}d{{z}
}-4\,{\it \Pi^2_{23}}\,vd{{y}}d{{z}}-$$$$-4\,{\it \Pi^1_{23}}\,ud{{y}}d{{z}}+2\,d{
{u}}d{{x}}+2\,d{{v}}d{{y}}+2\,d{{w}}d{{z}}.
\end{equation}

    The metric (\ref{6dflat}) is a flat $R_{i j k l}=0$ on solution of the KdF-equation (\ref{kdf2}).

  Simple modification of the metric $ds^2=(\ref{6dflat})$ leads to the following statement
\begin{theorem}
  The metric
$$
     \tilde ds^2=ds^2+\,\epsilon\,{{\it dy}}^{2}
$$
is the Richi -flat $(R_{ik}=0)$ on solution of the KdF-equation (\ref{kdf2}) but it has the nonzero curvature $R_{ijkl} \neq 0$.
\end{theorem}
\centerline{\bf References:}
\smallskip
\noindent 1. I.M. Krichever and S.P.Novikov, {\it Holomorphic bundles over algebraic curves and nonlinear equations}, Uspekhi Math. Nauk., v.35, 1980, p. 47--68.
\smallskip
\noindent 2. V.S. Dryuma, {\it K teorii  prostranstv postoyannoi krivizny}. Theor.and Math. Phys.,2006, v.146, N.1, p.42--54.
\noindent 3. N. Euler and M. Euler, {\it Multipotentialisation and iterating-solution formulae: The Krichever -Novikov equation }, Journal of Nonlinear mathematical Physics, V. 16, Suppl. (2009), 93--106.
\smallskip
\noindent 4. V. Dryuma, {\it The Riemann and Einsten-Weyl
geometries in theory of differential equations, their applications
and all that}. A.B.Shabat et all.(eds.), New Trends in
Integrability and Partial Solvability, Kluwer Academic Publishers,
Printed in the Netherlands , 2004, p.115--156.
\end{document}